\documentclass[aop,MSNbibl,citesort]{arximspdf}

\doi{10.1214/11-AOP648}
\volume{40}
\issue{3}
\pubyear{2012}
\firstpage{1285}
\lastpage{1315}

\makeatletter
\newcommand\R{{\mathbb{R}}}
\newcommand\C{{\mathbb{C}}}
\renewcommand\P{{\mathbf{P}}}
\newcommand\E{{\mathbf{E}}}
\newcommand\M{{\mathbf{M}}}
\renewcommand\Im{{\operatorname{Im}}}
\renewcommand\Re{{\operatorname{Re}}}
\newcommand\eps{{\varepsilon}}
\newcommand\dist{\operatorname{dist}}

 \newtheorem{theorem}{Theorem}
 \newtheorem{proposition}[theorem]{Proposition}
 \newtheorem{lemma}[theorem]{Lemma}
 \newtheorem{corollary}[theorem]{Corollary}
 \newproclaim{definition}[theorem]{Definition}
\newproclaim{remark}[theorem]{Remark}
\newcommand{\eqref}[1]{(\ref{#1})}
\makeatother

\begin{document}
\begin{frontmatter}

\title{Random covariance matrices: Universality of local statistics of
eigenvalues}
\runtitle{Universality for covariance matrices}

\begin{aug}
\author[A]{\fnms{Terence} \snm{Tao}\corref{}\thanksref{t1}\ead[label=e1]{tao@math.ucla.edu}}
\and
\author[B]{\fnms{Van} \snm{Vu}\thanksref{t2}\ead[label=e2]{vanvu@math.rutgers.edu}}
\runauthor{T. Tao and V. Vu}
\affiliation{UCLA and Rutgers University}
\address[A]{Department of Mathematics\\ UCLA\\
 Los Angeles, California
90095-1555\\
USA\\
\printead{e1}} 
\address[B]{Department of Mathematics\\
Rutgers University\\
 Piscataway, New Jersey 08854\\
USA\\
\printead{e2}}
\end{aug}
\thankstext{t1}{Supported by a grant from the MacArthur Foundation, by
NSF Grant DMS-06-49473, and by the NSF Waterman award.}
\thankstext{t2}{Supported by research Grants DMS-09-01216 and
AFOSAR-FA-9550-09-1-0167.}

\received{\smonth{12} \syear{2009}}
\revised{\smonth{1} \syear{2011}}

%
\begin{abstract}
We study the eigenvalues of the covariance matrix $\frac{1}{n} M^\ast
M$ of a~large rectangular matrix $M = M_{n,p} = (\zeta_{ij})_{1 \leq i
\leq p; 1 \leq j \leq n}$ whose entries are i.i.d. random variables of
mean zero, variance one, and having finite $C_0$th
moment for some sufficiently large constant $C_0$.

The main result of this paper is a Four Moment theorem for i.i.d.
covariance matrices (analogous to the Four Moment theorem for Wigner
matrices established by  the authors in  [\textit{Acta Math.} (2011) Random matrices: Universality of local eigenvalue statistics]
(see also [\textit{Comm. Math. Phys.}
\textbf{298} (2010) 549--572])). We can use this theorem together with existing results to
establish universality of local statistics of eigenvalues under mild conditions.

As a byproduct of our arguments, we also extend our previous results on
random Hermitian matrices to the case in which the entries have finite
$C_0$th moment rather than exponential decay.
\end{abstract}

%

\vspace*{-3pt}
\end{frontmatter}

\setcounter{footnote}{2}
\section{Introduction}
\vspace*{-3pt}

\subsection{The model}
The main purpose of this paper is to study the asymptotic local
eigenvalue statistics of covariance matrices of large random matrices.
Let us first fix the matrix ensembles that we will be studying.

\begin{definition}[(Random covariance matrices)] Let $n$ be a large
integer parameter going off to infinity, and let $p = p(n)$ be another
integer parameter such that $p \leq n$ and $\lim_{n \to\infty} p/n =
y$ for some $0 < y \leq1$. We let $M = M_{n,p} = (\zeta_{ij})_{1 \leq
i \leq p, 1 \leq j \leq n}$ be a random $p \times n$ matrix, whose
distribution is allowed to depend on $n$. We say that the matrix
ensemble $M$ \textit{obeys condition}  C1  with some exponent $C_0
\geq2$ if the random variables $\zeta_{ij}$ are jointly independent,
have mean\vadjust{\goodbreak} zero and variance $1$, and obey the moment condition $\sup
_{i,j} \E|\zeta_{ij}|^{C_0} \leq C$ for some constant $C$ independent
of $n,p$. We say that the matrix $M$ is \textit{i.i.d.} if the $\zeta
_{ij}$ are identically and independently distributed with law
independent of $n, p$.

Given such a matrix, we form the $n \times n$ \textit{covariance matrix}
$W = W_{n,p} := \frac{1}{n} M^\ast M$. This matrix has rank $p$ and so
the first $n-p$ eigenvalues are trivial; we order the (necessarily
positive) remaining eigenvalues of these matrices (counting
multiplicity) as
\[
0 \leq\lambda_1(W) \leq\cdots \leq\lambda_p(W).
\]
We often abbreviate $\lambda_i(W)$ as $\lambda_i$.
\end{definition}

Note that the only distributional hypothesis we require on the entries
$\zeta_{ij}$, besides the crucial joint independence hypothesis, are
moment conditions. In particular, we make no distinction between
continuous and discrete distributions here.

\begin{remark} In this paper, we will focus primarily on the case
$y=1$, but several of our results extend to other values of $y$ as
well. The case $p > n$ can be easily deduced from the $p<n$ case after
some minor notational changes by transposing the matrix $M$, which does
not affect the nontrivial eigenvalues of the covariance matrix. One can
also easily normalise the variance of the entries to be some other
quantity $\sigma^2$ than $1$ if one wishes. Observe that the
quantities $\sigma_i := \sqrt{n} \lambda_i^{1/2}$ can be interpreted
as the nontrivial singular values of the original matrix $M$, and
$\lambda_1,\ldots,\lambda_p$ can also be interpreted as the
eigenvalues of the $p \times p$ matrix $\frac{1}{n} M M^\ast$. It
will be convenient to exploit all three of these spectral
interpretations of $\lambda_1,\ldots,\lambda_p$ in this paper.
condition~{C1} is analogous to condition  {C0} for Wigner-type
matrices in \cite{TVlocal1}, but with the exponential decay hypothesis
relaxed to polynomial decay only.
\end{remark}

The well-known \textit{Marchenko--Pastur law} governs the bulk
distribution of the eigenvalues $\lambda_1,\ldots,\lambda_p$ of $W$:

\begin{theorem}[(Marchenko--Pastur law)] Assume condition  \textup{C1} with
\mbox{$C_0 > 2$}, and suppose that $p/n \to y$ for some $0 < y \leq1$. Then
for any $x > 0$, the random variables
\[
\frac{1}{p} | \{ 1 \leq i \leq p\dvtx  \lambda_i(W) \leq x \} |
\]
converge in probability to $\int_0^x \rho_{\mathrm{MP},y}(x)\,dx$, where
%
\begin{equation}\label{rhomp}
\rho_{\mathrm{MP},y}(x) := \frac{1}{2\pi xy} \sqrt{(b-x)(x-a)} 1_{[a,b]}(x)
\end{equation}
and
%
\begin{equation}\label{ab-def}
a := \bigl(1-\sqrt{y}\bigr)^2; \qquad b = \bigl(1+\sqrt{y}\bigr)^2.
\end{equation}
When furthermore $M$ is i.i.d., one can also obtain the case
$C_0=2$.\vadjust{\goodbreak}
\end{theorem}

\begin{pf} For the case $C_0 \geq4$, see \cite{marchenko,pastur}; for the case $C_0 > 2$, see \cite{wachter}; for the $C_0=2$
i.i.d. case, see \cite{yin}. Further results are known on the rate of
convergence: see \cite{gotze}.
\end{pf}

In this paper, we are concerned instead with the local eigenvalue
statistics. A~model case is the (\textit{complex}) \textit{Wishart ensemble}, in
which the $\zeta_{ij}$ are i.i.d. variables which are complex
Gaussians with mean zero and variance $1$. In this case, the
distribution of the eigenvalues $(\lambda_1,\ldots,\lambda_n)$ of
$W$ can be explicitly computed (as a special case of the \textit
{Laguerre unitary ensemble}). For instance, when $p=n$, the joint
distribution is given by the density function
%
\begin{equation} \label{jointcomplex}
\rho_n(\lambda_1,\ldots,\lambda_n) = c(n) \prod_{1 \le i < j \le
n} |\lambda_i -\lambda_j|^2 \exp\Biggl(-n\sum_{i=1}^n \lambda_i \Biggr)
\end{equation}
for some explicit normalization constant $c(n)$ whose exact value is
not important for this discussion.

Very similarly to the GUE case, one can use this explicit formula to
directly compute several local statistics,
including the distribution of the largest and smallest eigenvalues
\cite{edelman}, the correlation functions
\cite{nw} etc. Also in similarity to the GUE case, it is widely
conjectured that these statistics hold for a much larger class of
random matrices.
For some earlier results in this direction, we refer to \cite{soshnikov,TVhard,BenP,FS} and the references therein.

The goal of this paper is to establish a Four Moment theorem for random
covariance matrices, as an analogue of a recent result in \cite{TVlocal1}.
Roughly speaking, this theorem asserts that the asymptotic behaviour of
local statistics of the eigenvalues of $W_n$ are determined by the
first four moments of the entries.

\subsection{The Four Moment theorem}


To state the Four Moment theorem, we first need a definition.

\begin{definition}[(Matching)] We say that two complex random
variables~$\zeta$, $\zeta'$ \textit{match to order $k$} for some integer $k \geq
1$ if one has $\E\Re(\zeta)^m \Im(\zeta)^l = \E\Re(\zeta')^m
\Im(\zeta')^l$ for all $m,l \geq0$ with $m+l \leq k$.
\end{definition}

Our main result is the following.

\begin{theorem}[(Four Moment theorem)]\label{four-main1}
For sufficiently small $c_0>0$ and sufficiently large $C_0>0$
($C_0=10^4$ would suffice) the following holds for every $0 < \eps< 1$
and $k \geq1$.
Let $M = (\zeta_{ij})_{1 \leq i \leq p, 1 \leq j \leq n}$ and $M' =
(\zeta'_{ij})_{1 \leq i \leq p, 1 \leq j \leq n}$ be matrix ensembles
obeying condition \textup{C1} with the
the indicated constant $C_0$, and assume that for each $i,j$ that
$\zeta_{ij}$ and $\zeta'_{ij}$ match to order $4$. Let $W, W'$ be the
associated covariance matrices.
Assume also that $p/n \to y$ for some $0 < y
\leq1$.\vadjust{\goodbreak}

Let $G\dvtx  \R^k \to\R$ be a smooth function obeying the derivative bounds
%
\begin{equation}\label{G-deriv}
|\nabla^j G(x)| \leq n^{c_0}
\end{equation}
for all $0 \leq j \leq5$ and $x \in\R^k$.

Then for any $\eps p \le i_1 < i_2 <\cdots< i_k \le(1-\eps) p$, and
for $n$ sufficiently large depending on $\eps, k, c_0$ we have
%
\begin{equation} \label{eqnapproximation0}
 \quad |\E( G(n\lambda_{i_1}(W), \ldots , n\lambda_{i_k}(W))) -
\E( G(n\lambda_{i_1}(W'), \ldots , n\lambda_{i_k}(W')))| \le n^{-c_0}.
\end{equation}

If $\zeta_{ij}$ and $\zeta'_{ij}$ only match to order $3$ rather than $4$,
the conclusion \eqref{eqnapproximation0} still holds provided that
one strengthens \eqref{G-deriv} to
\[
|\nabla^j G(x)| \leq n^{-j c_1}
\]
for all $0 \leq j \leq5$ and $x \in\R^k$ and any $c_1 > 0$, provided
that $c_0$ is sufficiently small depending on $c_1$.
\end{theorem}

This is an analogue of \cite{TVlocal1}, Theorem 15, for covariance
matrices, with the main difference being that the exponential decay
condition from \cite{TVlocal1}, Theorem 15, has been weakened to the
high moment condition in  {C1}. This is achieved by an ``exponential
decay removing trick'' that relies on using a~truncated version of the
four moment theorem to extend the range of validity of a key ``gap
condition'' that is used in the proof of the above theorem. The same
trick also allows one to obtain a similar strengthening of the main
results of \cite{TVlocal1,TVlocal2}, thus relaxing the exponential
decay hypotheses in those results to high moment conditions. The value
$C_0=10^4$ is
\textit{ad hoc}, and we make no attempt to optimize this constant.

\begin{remark} \label{remarkoutsidethebulk} The reason that we
restrict the eigenvalues
to the bulk of the spectrum [$\eps p \leq i \leq(1-\eps) p$] is to
guarantee that the density function~$\rho_{\mathrm{MP},y}$ is bounded away from zero.
In view of the results in \cite{TVlocal2}, we expect that the result
extends to the edge of the spectrum as well. In particular, in view of
the results in \cite{BenP}, it is likely that the hard edge
asymptotics of Forrester~\cite{forrester} can be extended to a wider
class of ensembles. We will pursue this issue elsewhere.
\end{remark}

\begin{remark} As observed in \cite{ERSTVY}, the requirement that the
moments of $\zeta_{ij}$ and $\zeta'_{ij}$ match exactly can be
relaxed slightly. Indeed, to obtain the desired conclusions, it
suffices to require that for $k=1,2,3,4$, the $k${th} moments
of~$\zeta_{ij}$ and $\zeta'_{ij}$ differ by $O( n^{-(4-k)/2 - \delta}
)$ for some $\delta> 0$ independent of $n$. Indeed, if one inspects
the proof of the four moment theorem, and specifically the step in
which one performs a Taylor expansion argument to understand the effect
of exchanging a single entry $\zeta_{ij}$ with $\zeta'_{ij}$ on the
expectations in~\eqref{eqnapproximation0} (see~\cite{TVlocal1}, Section 3.2), the above near-matching property is sufficient to
ensure that this effect has magnitude $O( n^{-2-c} )$ for some $c>0$,
and so the net effect on \eqref{eqnapproximation0} after performing
$O(n^2)$ such exchange operations is acceptable. We omit the details.
This relaxed version of the four moment theorem is particularly useful
for dealing with Bernoulli distributions, which are completely
determined by their first four moments; see \cite{ERSTVY} for further
discussion.
\end{remark}

\subsection{Applications}

One can apply Theorem \ref{four-main1} in a similar way as its
counterpart \cite{TVlocal1}, Theorem 15, in order to obtain
universality results for large classes of random matrices. In many
cases, one can combine this theorem with existing partial results for
special ensembles to remove some of the moment assumptions. Let us
demonstrate this through an example concerning the universality of the
\textit{sine kernel}.

Using the explicit formula \eqref{jointcomplex}, Nagao and Wadati
\cite{nw} established the following result for the complex Wishart
ensemble, which roughly speaking asserts that the spectrum of such an
ensemble enjoys sine kernel statistics in the neighborhood of any bulk
energy level $0 < u < 4$.

\begin{theorem}[(Sine kernel for Wishart ensemble)]\label{wishart}
\cite{nw} Let $k \geq1$ be an integer, let $f\dvtx  \R^k \to\C$ be a
continuous function with compact support and symmetric with respect to
permutations, and let $0 < u < 4$; we assume all these quantities are
independent of $n$. Assume that\footnote{See Section \ref{notation-sec}
for the asymptotic notation we will be using.} $p = n +
O(1)$ (thus \mbox{$y=1$}), and that $W$ is given by the complex Wishart
ensemble. Let $\lambda_1,\ldots,\lambda_p$ be the nontrivial
eigenvalues of~$W$. Then the quantity
%
\begin{equation}\label{close}
\E\sum_{1 \leq i_1,\ldots,i_k \leq p, \hbox{ distinct}} f\bigl( n \rho_{\mathrm{MP},1}(u) (\lambda
_{i_1}-u),\ldots, n \rho_{\mathrm{MP},1}(u) (\lambda_{i_k}-u)\bigr)
\end{equation}
converges as $n \to\infty$ to
\[
\int_{\R^k} f(t_1,\ldots,t_k) \det( K(t_i,t_j) )_{1 \leq i,j \leq
k}\,dt_1 \cdots dt_k,
\]
where $K(x,y) := \frac{\sin(\pi(x-y))}{\pi(x-y)}$ is the sine kernel.
\end{theorem}

\begin{remark} The results in \cite{nw} allowed $f$ to be bounded
measurable rather than continuous, but when we consider discrete
ensembles later, it will be important to keep $f$ continuous.
\end{remark}

Returning to the bulk, the following extension was established by Ben
Arous and Pech\'{e} \cite{BenP}, as a variant of Johansson's result
\cite{Joh1} for random hermitian matrices. We say that a complex
random variable $\zeta$ of mean zero and variance one is \textit{Gauss
divisible} if $\zeta$ has the same distribution as
$\zeta= (1-t)^{1/2} \zeta' + t^{1/2} \zeta''$ for some $0 < t < 1$
and some independent random variables~$\zeta'$, $\zeta''$ of mean zero
and variance $1$, with $\zeta''$ distributed according to the complex
Gaussian.\vadjust{\goodbreak}

\begin{theorem}[(Sine kernel for Gaussian divisible ensemble)]\label{gauss-divisible}
\cite{BenP} Theorem~\ref{wishart} [which is for the
Wishart ensemble and for $p=n+O(1)$] can be extended to the case when
$p = n + O(n^{43/48})$ (so $y$ is still $1$), and when $M$ is an i.i.d.
matrix obeying condition \textup{C1} with $C_0=2$, and with the $\zeta
_{ij}$ gauss divisible.
\end{theorem}

Using Theorem \ref{four-main1} and Theorem \ref{gauss-divisible}
(in exactly the same way we used \cite{TVlocal1}, Theorem 15, and
Johansson's theorem \cite{Joh1} to establish \cite{TVlocal1}, Theorem 11), we can extend Theorem \ref{gauss-divisible} from the
gauss divisible case to a more general situation.

\begin{corollary}[(Sine kernel for more general ensembles)]\label{main-thm} Theorem \ref{wishart} can be extended to the case when $p =
n + O(n^{43/48})$ (so $y$ is still $1$), and when $M$ is an i.i.d.
matrix obeying condition \textup{C1} with $C_0$ sufficiently large ($C_0 =
10^4$ would suffice), and where the real and imaginary parts of $\zeta
_{ij}$ are i.i.d. and are supported on at least three points.
\end{corollary}

\begin{pf} (Sketch) It was shown in \cite{TVlocal1}, Corollary 30,
that if the real and imaginary parts of a complex random variable
$\zeta$ were independent with mean zero and variance one, and both
were supported on at least three points, then $\zeta$ matched to order
$4$ with a gauss divisible random variable~$\zeta'$ with finite $C_0$
moment (indeed, if one inspects the convexity argument used to solve
the moment problem in \cite{TVlocal1}, Lemma 28, the Gauss divisible
random variable could be taken to be the sum of a Gaussian variable and
a discrete variable, and in particular is thus exponentially decaying).
If one lets $M'$ be the i.i.d. matrix whose coefficients have entries
$\zeta'$, then Theorem~\ref{gauss-divisible} asserts that the
conclusions of Theorem \ref{wishart} hold for $M'$. Using Theorem \ref{four-main1} exactly as in the proof of \cite{TVlocal1}, Theorem 11,
(and approximating $f$ uniformly by smooth functions), we conclude that
the conclusions of Theorem \ref{wishart} hold for $M$ also.
\end{pf}

One can also extend the above argument to cover cases in which the real
and imaginary parts of $\zeta_{ij}$ are not i.i.d. by an analysis of
the moment matching problem for complex random variables (and in
particular, by extending the three-moment analysis in Lemma \ref{lemmamatching} below to four moments), but we will not do so here.

The arguments in this paper will be a nonsymmetric version of those in
\cite{TVlocal1}. The arguments in \cite{TVlocal1} started with
analyzing the stability of the eigenvalue equation
$M v_i =\lambda_i v _i$ where $M$ is a random Hermitian matrix and
$\lambda_i$ is the $i$th eigenvalue with
eigenvector $v$. For the situation considered in this paper, it is
tempting to similarly analyze the eigenvalue equation $Wv_i =\lambda_i
v_i$ for the covariance matrix~$W$. However, this does not work, since
the covariance\vadjust{\goodbreak} matrix $W$, while random, does not have independent
entries. The new idea here is to work with a system of two equations
%
\begin{equation}\label{mu}
M u_i = \sigma_i v_i
\end{equation}
and
%
\begin{equation}\label{muv}
M^{\ast} v_i= \sigma_i u_i,
\end{equation}
where $u_i$ and $v_i$ are the left and right singular vectors of $M$.
This leads to a~number of technical issues that need to be addressed
through the paper.

One can combine the singular value equations \eqref{mu}, \eqref{muv}
into a single eigenvalue equation
\[
\M
\pmatrix{\displaystyle v_i \cr\displaystyle  u_i
}
= \sigma_i
\pmatrix{\displaystyle v_i \cr\displaystyle  u_i
},
\]
where $\M$ is the \textit{augmented matrix}
%
\begin{equation}\label{augment}
\M:=
\pmatrix{\displaystyle 0 & M \cr\displaystyle  M^\ast& 0
}
.
\end{equation}
Thus one can view the singular values of an i.i.d. matrix as being
essentially given by the eigenvalues of a slightly larger Hermitian
matrix which is of Wigner type except that the entries have been zeroed
out on two diagonal blocks. We will take advantage of thus augmented
perspective in some parts of the paper (particularly when we wish to
import results from \cite{TVlocal1} as black boxes), but in other
parts it will in fact be more convenient to work with $M$ directly. In
particular, the fact that many of the entries in \eqref{augment} are
zero (and in particular, have zero mean and variance) seems to make it
difficult to directly apply parts of the arguments from \cite{TVlocal1} (particularly those that are probabilistic in
nature,\footnote{A typical instance of a probabilistic argument that
encounters difficulty when there are many zero entries arises when one
wants to estimate the distance $\dist(X,V)$ between a~random vector $X
= (\xi_1,\ldots,\xi_n)$ (which one should think of as something like
a row of~$\M$) and a fixed subspace $V$. If all the entries of $X$ are
i.i.d. with mean zero and constant variance, then an easy second moment
computation allows one to control $\E\dist(X,V)^2$ exactly in terms
of the codimension of $V$; in particular, no knowledge of the
orientation of $V$ is required. One also obtains reasonable upper and
lower bounds on this quantity if the variance is not constant, but is
also bounded above and below. However, if many of the entries of $X$
have zero variance (i.e., they vanish), then one has difficulty lower
bounding~$\E\dist(X,V)^2$ because one has to somehow exclude the
possibility that the normal vectors to $V$ have almost all of their
$\ell^2$ mass supported on those zero variance entries. We do not know
how to address this problem in general. \textit{Note added in proof}:
Several months after the submission of this paper, Erd\H{o}s, Yau  and
Yin \cite{EYY0,EYY} were able to obtain universality results for some
classes of generalized Wigner matrices (such as band matrices) in which
some entries are permitted to have zero variance. However, one of their
key assumptions is that the matrix of (normalised) variances has a
simple eigenvalue at~$1$, and this assumption does not hold for the
augmented matrix \eqref{augment}.} rather than deterministic)
directly to the augmented matrix,\vadjust{\goodbreak} and will instead work with $M$
directly in these cases. Nevertheless, one can view this connection as
a heuristic explanation as to why some (but not all) of the machinery
in the Hermitian eigenvalue problem can be transferred to the
non-Hermitian singular value problem.

\subsection{Extensions}

In a very recent work, Erd\H{o}s et al. \cite{ESYY}
extended\footnote{Even more recently, a similar result was also
established by P\'{e}ch\'{e} \cite{peche}.} Theorem \ref{gauss-divisible} to a large class of matrices, assuming that the
distribution of the entries $\zeta_{ij}$ is sufficiently smooth and
obeys a log-Sobolev inequality. While their results do not apply for
entries with discrete distributions, it allows one to extend Theorem
\ref{gauss-divisible} to the case when $t$ is a negative power of $n$.
Given this, one can use the argument in \cite{ERSTVY} to
remove the requirement that the real and imaginary parts of $\zeta
_{ij}$ be supported on at least three points.

We can also have the following analogue of \cite{ERSTVY}, Theorem 2.

\begin{theorem}[(Universality of averaged correlation function)]\label{main1}
Fix $\eps>0$ and $u$ such that
$0 < u - \eps< u+\eps< 4$. Let $k \geq1$ and
let $f\dvtx  \R^k \to\R$ be a~continuous, compactly supported function,
and let $W=W_{n,n}$ be a random covariance matrix, with $n$ assumed
large depending on $u,\eps,k$. Then the quantity
%
\begin{eqnarray}\label{episode}
&&\frac{1}{2\eps} \int_{u-\eps}^{u+\eps} \int_{\R^k} f(t_1,\ldots,t_k)
\frac{1}{(n\rho_{\mathrm{MP},1}(u'))^k} p_n^{(k)}\nonumber\hspace*{-25pt}
\\[-8pt]
\\[-8pt]
&&\hphantom{\frac{1}{2\eps} \int_{u-\eps}^{u+\eps} \int_{\R^k} }{}\times\biggl(u'+\frac{t_1}{n \rho
_{\mathrm{MP},1}(u')},
\ldots, u'+\frac{t_k}{n \rho_{\mathrm{MP},1}(u')} \biggr)\,dt_1 \cdots dt_k \,du'
\nonumber\hspace*{-25pt}
\end{eqnarray}
converges as $n \to\infty$ to
\[
\int_{\R^k} f(t_1,\ldots,t_k)
\det( K(t_i,t_j) )_{i,j=1}^k\,dt_1 \cdots dt_k,
\]
where $K(x,y)$ is the \textit{Dyson sine kernel}
%
\begin{equation}\label{dyson}
K(x,y) := \frac{\sin( \pi(x-y) )}{\pi(x-y)},
\end{equation}
and the $k$-point correlation function $p_n^{(k)}\dvtx  \R^k \to\R^+$ is
the unique symmetric probability distribution such that
\[
\int_{\R^k} f(\alpha_1,\ldots,\alpha_k) p_n^{(k)}(\alpha_1,\ldots
,\alpha_k) := k! \sum_{1 \leq i_1 < \cdots < i_k \leq n} f( \lambda
_1, \ldots, \lambda_n )
\]
for all symmetric test functions $f$. (If $W$ is a discrete ensemble,
one has to interpret $p_n^{(k)}$ as a distribution or a probability
measure rather than as a~function.)
\end{theorem}

The detailed proof of Theorem \ref{main1} are essentially the same as
the proof of~\cite{ERSTVY}, Theorem 2, and is omitted.\vadjust{\goodbreak}

\begin{remark} The four moment theorem controls the distribution of
individual eigenvalues (or singular values) $\lambda_i$ , but as
indicated above, this control can then be used to obtain control of
correlation expressions such as~\eqref{episode}. The local relaxation
flow methods developed in \cite{ESY1,ESY2,ESY3,ERSY,ERSY2,ESYY}, by contrast, are
focused on individual energy levels $u$ rather than individual
eigenvalues. As such, they provide an alternate approach to controlling
correlation expressions such as \eqref{episode}, but we do not know
how to convert such information back to control on individual
eigenvalues or singular values in general, because the standard
deviation of each eigenvalue can exceed (by a logarithmic factor, see
\cite{Gus}) the scale of the mean eigenvalue spacing, which is the
scale at which the correlation estimates operate at.
\end{remark}

\subsection{Notation}\label{notation-sec}

Throughout this paper, $n$ will be an asymptotic parameter going to
infinity. Some quantities (e.g., $\eps$, $y$ and $C_0$) will remain
independent of $n$, while other quantities (e.g., $p$, or the matrix
$M$) will depend on $n$.
All statements here are understood to hold only in the asymptotic
regime when $n$ is sufficiently large depending on all quantities that
are independent of $n$. We write $X = O(Y)$, $Y = \Omega(|X|)$, $|X|
\ll Y$, or $Y \gg|X|$ if one has $|X| \leq CY$ for all sufficiently
large $n$ and some $C$ independent of $n$. [Note however that $C$ is
allowed to depend on other quantities independent of $n$, such as $\eps
$ and $y$, unless otherwise stated; we will sometimes emphasise this by
using subscripts, thus, for instance, $X = O_a(Y)$ denotes the estimate
$|X| \leq C_a Y$ for some constant $C$ depending only on $a$.] We write
$X = o(Y)$ if $|X| \leq c(n) Y$ where $c(n) \to0$ as $n \to\infty$.
We write $X = \Theta(Y)$ if $X \ll Y \ll X$, thus, for instance, if $p/n
\to y$ for some $0 < y \leq1$ then $p = \Theta(n)$.

We write $\sqrt{-1}$ for the complex imaginary unit, in order to free
up the letter $i$ to denote an integer (usually between $1$ and $n$).

We write $\|X\|$ for the length of a vector $X$, $\| A \| = \| A \|
_{op}$ for the operator norm of a matrix $A$, and $\|A\|_F =
\operatorname{tr}(AA^*)^{1/2}$ for the Frobenius (or Hilbert--Schmidt) norm.

We will need to quantify the intuitive assertion that a given event $E$
occurs ``frequently,'' as follows.

\begin{definition}[(Frequent events)]\label{freq-def}
$\!\!\!$\cite{TVlocal1}
Let $E$ be an event depending~on~$n$.
\begin{itemize}
\item
$E$ holds \textit{with high probability} if $\P(E) \geq
1-O(n^{-c})$ for some constant $c>0$ (independent of $n$).
\item
$E$ holds \textit{with overwhelming probability} if $\P(E) \geq
1-O_C(n^{-C})$ for \textit{every} constant $C>0$.
\item
$E$ holds \textit{almost surely} if $\P(E)=1$.
\end{itemize}
\end{definition}

\section{The gap property and the exponential decay removing trick}

The following property, which roughly speaking asserts that
unexpectedly small eigenvalue spacings are rare, plays an important
role in proving the main results of \cite{TVlocal1}.\vadjust{\goodbreak}

\begin{definition}[(Gap property)] Let $M$ be a matrix ensemble obeying
condition  {C1}. We say that $M$ obeys the \textit{gap property} if
for every $\eps, c > 0$ (independent of $n$), and for every $\eps p
\leq i \leq(1-\eps) p$, one has $|\lambda_{i+1}(W) - \lambda_i(W)|
\geq n^{-1-c}$ with high probability. (The implied constants in this
statement are allowed to depend on $\eps$ and $c$.)
\end{definition}

In the Wigner case, it was shown that exponential decay of the atom
distribution implied the gap property, and the gap property was then
used to establish deduce the four moment theorem from a ``truncated
four moment theorem.'' As it turns out, the proof of this latter
theorem does not require exponential decay of the atom distribution,
relying instead on the weaker hypothesis that a sufficiently high
moment of the atom distribution is finite. A new technical observation
of this paper is that one can use the truncated four moment theorem to
extend the gap property from exponentially decaying atom distributions
to distributions with sufficiently high moments finite, and as a
consequence we can extend the full Four Moment theorem to this case also.

We turn to the details. First, as an analogue of \cite{TVlocal1}, Theorem 19, we prove the following theorem, using a slight
modification of the method in \cite{TVlocal1}.

\begin{theorem}[(Gap theorem)]\label{gap}
$\!\!\!$Let $M = (\zeta_{ij})_{1 \leq
i \leq p, 1 \leq j \leq n}$ obey condition~\textup{C1} for some $C_0$, and
suppose that the coefficients $\zeta_{ij}$
are \textit{exponentially decaying} in the sense that $\P( |\zeta
_{ij}| \geq t^C ) \leq\exp(-t)$ for all $t \geq C'$ for all $i,j$ and
some constants~$C$, $C'>0$. Then $M$ obeys the gap property.
\end{theorem}

Next, we have the following analogue of \cite{TVlocal1}, Theorem 15.

\begin{theorem}[(Four Moment theorem with Gap assumption)]\label{four-main}
For sufficiently small $c_0>0$ and sufficiently large $C_0>0$
($C_0=10^4$ would suffice) the following holds for every $0 < \eps< 1$
and $k \geq1$.
Let $M = (\zeta_{ij})_{1 \leq i \leq p, 1 \leq j \leq n}$ and $M' =
(\zeta'_{ij})_{1 \leq i \leq p, 1 \leq j \leq n}$ be matrix ensembles
obeying condition \textup{C1} with
the indicated constant $C_0$, and assume that for each $i,j$ that
$\zeta_{ij}$ and $\zeta'_{ij}$ match to order~$4$. Let $W, W'$ be the
associated covariance matrices. Assume also that~$M$ and
$M'$ obeys the gap property, and that $p/n \to y$ for some $0 < y \leq1$.

Let $G\dvtx  \R^k \to\R$ be a smooth function obeying the derivative bounds
%
\begin{equation}\label{G-deriv-0}
|\nabla^j G(x)| \leq n^{c_0}
\end{equation}
for all $0 \leq j \leq5$ and $x \in\R^k$.

Then for any $\eps p \le i_1 < i_2 <\cdots< i_k \le(1-\eps) p$, and
for $n$ sufficiently large depending on $\eps, k, c_0$ we have
%
\begin{equation} \label{eqnapproximation}
 |\E( G(n\lambda_{i_1}(W), \ldots , n\lambda_{i_k}(W))) -
\E( G(n\lambda_{i_1}(W'), \ldots , n\lambda_{i_k}(W')))| \le n^{-c_0}.\hspace*{-30pt}
\end{equation}

If $\zeta_{ij}$ and $\zeta'_{ij}$ only match to order $3$ rather than $4$,
the conclusion \eqref{eqnapproximation} still holds provided that one
strengthens \eqref{G-deriv-0} to
\[
|\nabla^j G(x)| \leq n^{-j c_1}
\]
for all $0 \leq j \leq5$ and $x \in\R^k$ and any $c_1 > 0$, provided
that $c_0$ is sufficiently small depending on $c_1$.
\end{theorem}

This theorem is weaker than Theorem \ref{four-main1}, as we assume the
gap property. Besides the fact that we consider singular values here
instead of eigenvalues, the main difference between this result and
\cite{TVlocal1}, Theorem 15, is that in the latter we assume
exponential decay rather than the gap property. However, this
difference is only a formality, since in the proof of \cite{TVlocal1}, Theorem 15, the only place we used exponential decay is to prove the
gap property (via \cite{TVlocal1}, Theorem 19).

The core of the proof of Theorem \ref{four-main} is a truncated four
moment theorem (Theorem \ref{four-main-truncate}), which allows us to insert
information such as the gap property into the test function $G$.

By combining Theorem \ref{four-main} with Theorem \ref{gap}, we
obtain Theorem \ref{four-main1} in the case when the coefficients
$\zeta_{ij}$ are exponentially decaying. To remove the exponential
decay hypothesis, we will apply the truncated four moment theorem
(Theorem~\ref{four-main-truncate}) a second time, together with a
moment matching argument (Lemma \ref{lemmamatching}) to eliminate
this hypothesis from Theorem \ref{gap}.

\begin{theorem}[(Gap theorem)]\label{gap1} Assume that $M = (\zeta
_{ij})_{1 \leq i \leq p, 1 \leq j \leq n}$ satisfies condition \textup{C1}
with $C_0$ sufficiently large.
Then $M$ obeys the gap property.
\end{theorem}

Theorem \ref{four-main1} follows directly from Theorems \ref{four-main} and \ref{gap1}.

The rest of the paper is organized as follows. The next three sections
are devoted to technical lemmas. The proofs of Theorems
\ref{four-main} and \ref{gap1} are presented in Section~\ref{fmt},
assuming Theorems \ref{four-main-truncate} and \ref{gap}.
The proofs of these latter two theorems
are presented in Sections \ref{sectionfmtr} and \ref{sectiongap},
respectively.


\section{The main technical lemmas}

\textit{Important note.}
The arguments in this paper are very similar
to, and draw heavily from, the previous paper \cite{TVlocal1} of the
authors. We recommend therefore that the reader be familiar with that
paper first, before reading the current one.

In the proof of the Four Moment theorem (as well as the Gap theorem)
for $n \times n$ Wigner matrices in \cite{TVlocal1}, a crucial
ingredient was a variant of the
\textit{Delocalization Theorem} of Erd\"{o}s, Schlein  and Yau \cite{ESY1,ESY2,ESY3}. This result asserts (assuming uniformly exponentially
decaying distribution for the coefficients) that with overwhelming
probability, all the unit eigenvectors of the Wigner matrix have
coefficients $O(n^{-1/2+o(1)})$ (thus, the ``$\ell^2$ energy'' of the
eigenvector is spread out more or less uniformly amongst the $n$
coefficients). When one just assumes uniformly bounded $C_0$ moment
rather than uniform exponential decay, the bound becomes $O(
n^{-1/2+O(1/C_0)} )$ instead (where the implied constant in the
exponent is uniform in $C_0$).

Similarly, to prove the Four Moment and Gap theorems in this paper, we
will need a Delocalization theorem for the \textit{singular vectors} of
the matrix~$M$. We define a \textit{right singular vector} $u_i$
(resp.,
\textit{left singular vector} $v_i$) with singular value $\sigma_i(M) =
\sqrt{n} \lambda_i(W)^{1/2}$ to be an eigenvector of $W = \frac
{1}{n} M^\ast M$ (resp., $\tilde W = \frac{1}{n} M M^\ast$) with
eigenvalue $\lambda_i$.
In the \textit{generic case} when the singular values are simple (i.e.,
$0 < \sigma_1 < \cdots < \sigma_p$), we observe from the singular
value decomposition that one can find orthonormal bases $u_1,\ldots
,u_p \in\C^n$ and $v_1,\ldots,v_p \in\C^p$ for the corange $\ker
(M)^\perp$ of $M$ and of $\C^p$, respectively, such that
\[
M u_i = \sigma_i v_i
\]
and
\[
M^\ast v_i = \sigma_i u_i.
\]
Furthermore, in the generic case the unit singular vectors $u_i, v_i$
are determined up to multiplication by a complex phase $e^{i\theta}$.

We will establish the following Erd\"{o}s--Schlein--Yau type
delocalization theorem (analogous to \cite{TVlocal1}, Proposition 62),
which is an essential ingredient to Theorems~\ref{four-main}, \ref{gap} and is also of some independent interest.

\begin{theorem}[(Delocalization theorem)]\label{delocal} Suppose that
$p/n \to y$ for some $0 < y \leq1$, and
let $M$ obey condition \textup{C1} for some $C_0 \geq2$. Suppose further
that that $|\zeta_{ij}| \leq K$ almost surely for some $K > 1$ (which
can depend on~$n$) and all $i,j$, and that the probability distribution
of $M$ is continuous.
Let $\eps> 0$ be independent of $n$. Then with overwhelming
probability, all the unit left and right singular vectors of $M$ with
eigenvalue $\lambda_i$ in the interval $[a+\eps,b-\eps]$ [with $a,b$
defined in \eqref{ab-def}] have all coefficients uniformly of size $O(
K n^{-1/2} \log^{10} n)$.
\end{theorem}

The factors $K \log^{10} n$ can probably be improved slightly, but
anything which is polynomial in $K$ and $\log n$ will suffice for our
purposes. Observe that if $M$ obeys condition  {C1}, then each
event $|\zeta_{ij}| \leq K$ with $K := n^{10/C_0}$ (say) occurs with
probability $1-O(n^{-10})$. Thus, in practice, we will be able to apply
the above theorem with $K = n^{10/C_0}$ without difficulty. The
continuity hypothesis is a technical one, imposed so that the singular
values are almost surely simple, but in practice we will be able to
eliminate this hypothesis by a limiting argument (as none of the bounds
will depend on any quantitative measure of this continuity).

As with other proofs of delocalization theorems in the literature,
Theorem~\ref{delocal} is in turn deduced from the following eigenvalue
concentration bound (analogous to \cite{TVlocal1}, Proposition 60).

\begin{theorem}[(Eigenvalue concentration theorem)]\label{eigenconc} Let
the hypotheses be as in Theorem \ref{delocal}, and let $\delta> 0$ be
independent of $n$. Then for any interval $I \subset[a+\eps, b-\eps
]$ of length $|I| \geq K^2 \log^{20} n / n$, one has with overwhelming
probability (uniformly in $I$) that
\[
\biggl|N_I - p \int_I \rho_{\mathrm{MP},y}(x)\,dx\biggr| \leq\delta p,
\]
where
%
\begin{equation}\label{NI-def}
N_I := \{ 1 \leq i \leq p\dvtx  \lambda_i(W) \in I \}
\end{equation}
is the number of eigenvalues in $I$.
\end{theorem}

We remark that a very similar result (with slightly different
hypotheses on the parameters and on the underlying random variable
distributions) was recently established in \cite{ESYY}, Corollary 7.2.

We isolate one particular consequence of Theorem \ref{eigenconc} (also
established in~\cite{GZ}):

\begin{corollary}[(Concentration of the bulk)]\label{bulk} Let the
hypotheses be as in Theorem \ref{delocal}. Then there exists $\eps' >
0$ independent of $n$ such that with overwhelming probability, one has
$a+\eps' \leq\lambda_i(W) \leq b-\eps'$ for all $\eps p \leq i \leq
(1-\eps) p$.
\end{corollary}

\begin{pf} From Theorem \ref{eigenconc}, we see with overwhelming
probability that the number of eigenvalues in $[a+\eps',b-\eps']$ is
at least $(1-\eps) p$, if $\eps'$ is sufficiently small depending on
$\eps$. The claim follows.
\end{pf}

\section{Basic tools}

\subsection{Tools from linear algebra}\label{tools-sec}

In this section, we recall some basic identities and inequalities from
linear algebra which will be used in this paper.

We begin with the Cauchy interlacing law and the Weyl inequalities.

\begin{lemma}[(Cauchy interlacing law)]\label{cauchy} Let $1 \leq p \leq n$.
\begin{longlist}[(iii)]
\item[(i)] If $A_n$ is an $n \times n$ Hermitian matrix, and $A_{n-1}$
is an $n-1 \times n-1$ minor, then $\lambda_i(A_n) \leq\lambda
_i(A_{n-1}) \leq\lambda_{i+1}(A_n)$ for all $1 \leq i < n$.
\item[(ii)] If $M_{n,p}$ is a $p \times n$ matrix, and $M_{n,p-1}$ is
an $p-1 \times n$ minor, then $\sigma_i(M_{n,p}) \leq\sigma_i(
M_{n,p-1} ) \leq\sigma_{i+1}(M_{n,p})$ for all $1 \leq i < p$.
\item[(iii)] If $p<n$, if $M_{n,p}$ is a $p \times n$ matrix, and
$M_{n-1,p}$ is a $p \times n-1$ minor, then $\sigma_{i-1}(M_{n,p})
\leq\sigma_i(M_{n-1,p}) \leq\sigma_i(M_{n,p})$ for all $1 \leq i
\leq p$, with the understanding that $\sigma_0(M_{n,p})=0$. [For
$p=n$, one can also use the transpose of~(\textup{ii}) instead.]
\end{longlist}
\end{lemma}

\begin{pf} Claim (i) follows from the minimax formula
\[
\lambda_i(A_n) = \inf_{V:  \dim(V)=i} \sup_{v \in V:  \|v\|=1} v^\ast
A_n v,
\]
where $V$ ranges over $i$-dimensional subspaces in $\C^n$. Similarly,
(ii) and (iii) follow from the minimax formula
\[
\sigma_i(M_{n,p}) = \inf_{V:  \dim(V)=i+n-p} \sup_{v \in V:  \|v\|=1}
\|M_{n,p} v\|.\vspace*{-2pt}
\]
\upqed
\end{pf}

\begin{lemma}[(Weyl inequality)]\label{weyl}
Let $1 \leq p \leq n$.
\begin{itemize}
\item If $A, B$ are $n \times n$ Hermitian matrices, then $\|\lambda
_i(A) - \lambda_i(B)| \leq\|A-B\|_{op}$ for all $1 \leq i \leq n$.
\item If $M, N$ are $p \times n$ matrices, then $\|\sigma_i(M) -
\sigma_i(N)| \leq\|M-N\|_{op}$ for all $1 \leq i \leq p$.\vspace*{-2pt}
\end{itemize}
\end{lemma}

\begin{pf} This follows from the same minimax formulae used to
establish Lemma \ref{cauchy}.\vspace*{-2pt}
\end{pf}

\begin{remark} One can also deduce the singular value versions of
Lemmas~\ref{cauchy}, \ref{weyl} from their Hermitian counterparts by
using the augmented matrices \eqref{augment}. We omit the details.\vspace*{-2pt}
\end{remark}

We have the following elementary formula for a component of an
eigenvector of a Hermitian matrix, in terms of the eigenvalues and
eigenvectors of a minor.\vspace*{-2pt}

\begin{lemma}[(Formula for coordinate of an eigenvector)]\label{lemmafirstcoordinate}
\cite{ESY1} Let
\[
A_n =
\pmatrix{\displaystyle A_{n-1} & X \cr\displaystyle  X^\ast& a
}
\]
be a $n \times n$ Hermitian matrix for some $a \in\R$ and $X \in\C
^{n-1}$, and let $
{ v \choose x
}
$ be a~unit eigenvector of $A_n$ with eigenvalue $\lambda_i(A_n)$,
where $x \in\C$ and $v \in\C^{n-1}$. Suppose that none of the
eigenvalues of $A_{n-1}$ are equal to $\lambda_i(A_n)$. Then
\[
|x|^2 = \frac{1}{1 + \sum_{j=1}^{n-1} (\lambda_j(A_{n-1})-\lambda
_i(A_n))^{-2} |u_j(A_{n-1})^\ast X|^2},
\]
where $u_1(A_{n-1}),\ldots,u_{n-1}(A_{n-1}) \in\C^{n-1}$ is an
orthonormal eigenbasis corresponding to the eigenvalues $\lambda
_1(A_{n-1}),\ldots,\lambda_{n-1}(A_{n-1})$ of $A_{n-1}$.\vspace*{-2pt}
\end{lemma}

\begin{pf}
See, for example, \cite{TVlocal1}, Lemma 41.\vspace*{-2pt}
\end{pf}

This implies an analogous formula for singular vectors.\vspace*{-2pt}

\begin{corollary}[(Formula for coordinate of a singular vector)]\label{singcoord}
Let \mbox{$p,n \geq1$}, and let
\[
M_{p,n} =
\pmatrix{\displaystyle M_{p,n-1} & X
}\vadjust{\goodbreak}
\]
be a $p \times n$ matrix for some $X \in\C^p$, and let $
{ u \choose x
}
$ be a right unit singular vector of $M_{p,n}$ with singular value
$\sigma_i(M_{p,n})$, where $x \in\C$ and $u \in\C^{n-1}$. Suppose
that none of the singular values of $M_{p,n-1}$ are equal to $\sigma
_i(M_{p,n})$. Then
\[
|x|^2 = \Biggl(1 + \sum_{j=1}^{\min(p,n-1)} \frac{\sigma
_j(M_{p,n-1})^2}{(\sigma_j(M_{p,n-1})^2-\sigma_i(M_{p,n})^2)^{2}}
|v_j(M_{p,n-1})^\ast X|^2\Biggr)^{-1},
\]
where $v_1(M_{p,n-1}), \ldots, v_{\min(p,n-1)}(M_{p,n-1}) \in\C^p$
is an orthonormal system of left singular vectors corresponding to the
nontrivial singular values of~$M_{p,n-1}$.

In a similar vein, if
\[
M_{p,n} =
\pmatrix{\displaystyle M_{p-1,n} \cr\displaystyle  Y^\ast
}
\]
for some $Y \in\C^n$, and $
\left({ v \enskip y
}\right)
$ is a left unit singular vector of $M_{p,n}$ with singular value
$\sigma_i(M_{p,n})$, where $y \in\C$ and $v \in\C^{p-1}$, and none
of the singular values of $M_{p-1,n}$ are equal to $\sigma
_i(M_{p,n})$, then
\[
|y|^2 = \Biggl(1 + \sum_{j=1}^{\min(p-1,n)} \frac{\sigma
_j(M_{p-1,n})^2}{(\sigma_j(M_{p-1,n})^2-\sigma_i(M_{p,n})^2)^{2}}
|u_j(M_{p-1,n})^\ast Y|^2\Biggr)^{-1},
\]
where $u_1(M_{p-1,n}),\ldots, u_{\min(p-1,n)}(M_{p-1,n}) \in\C^n$
is an orthonormal system of right singular vectors corresponding to the
nontrivial singular values of~$M_{p-1,n}$.
\end{corollary}

\begin{pf} We just prove the first claim, as the second is proven
analogously (or by taking adjoints). Observe that
$
{ u \choose x
}
$ is a unit eigenvector of the matrix
\[
M^\ast_{p,n} M_{p,n} =
\pmatrix{\displaystyle
M_{p,n-1}^\ast M_{p,n-1}
&
M_{p,n-1}^\ast X \cr\displaystyle  X^\ast M_{p,n-1} & |X|^2 }
\]
with eigenvalue $\sigma_i(M_{p,n})^2$. Applying Lemma \ref{lemmafirstcoordinate}, we obtain
\begin{eqnarray*}
|x|^2 &=& \Biggl(1 + \sum_{j=1}^{n-1} \bigl(\lambda_j(M_{p,n-1}^\ast
M_{p,n-1})-\sigma_i(M_{p,n})^2\bigr)^{-2}\\
&&\hspace*{7pt}\hphantom{\Biggl(1 + \sum_{j=1}^{n-1}}{}\times |u_j(M_{p,n-1}^\ast
M_{p,n-1})^\ast M_{p,n-1}^\ast X|^2\Biggr)^{-1}.
\end{eqnarray*}
But $u_j(M_{p,n-1}^\ast M_{p,n-1})^\ast M_{p,n-1}^\ast= \sigma
_j(M_{p,n-1}) v_j(M_{p,n-1})^\ast$ for the $\min(p,n-1)$ nontrivial
singular values (possibly after relabeling the $j$), and vanishes for
trivial ones, and $\lambda_j(M_{p,n-1}^\ast M_{p,n-1}) = \sigma
_j(M_{p,n-1})^2$, so the claim follows.
\end{pf}

The \textit{Stieltjes transform} $s(z)$ of a Hermitian matrix $W$ is
defined for complex $z$ by the formula
\[
s(z) := \frac{1}{n} \sum_{i=1}^n \frac{1}{\lambda_i(W)-z}.
\]
It has the following alternate representation (see, e.g., \cite{BS}, Chapter 11).

\begin{lemma}\label{stielt} Let $W = (\zeta_{ij})_{1 \leq i,j \leq
n}$ be a~Hermitian matrix, and let $z$ be a~complex number not in the
spectrum of $W$. Then we have
\[
s_n(z) = \frac{1}{n} \sum_{k=1}^n \frac{1}{\zeta_{kk} - z -
a_k^\ast(W_k - zI)^{-1} a_k},
\]
where $W_k$ is the $n-1 \times n-1$ matrix with the $k${th} row and column removed, and $a_k \in\C^{n-1}$ is the
$k$th column of $W$ with the $k${th} entry removed.
\end{lemma}

\begin{pf} By Schur's complement, $\frac{1}{\zeta_{kk} - z -
a_k^\ast(W_k - zI)^{-1} a_k}$ is the $k$th diagonal
entry of $(W-zI)^{-1}$. Taking traces, one obtains the claim.
\end{pf}

\subsection{Tools from probability theory}\label{prob-sec}

We will rely frequently on the following concentration of measure
result for projections of random vectors.

\begin{lemma}[(Distance between a random vector and a subspace)]\label{lemmaprojection} Let $X=(\xi_1, \ldots , \xi_n) \in\C^n$ be a
random vector whose entries are independent with mean zero, variance $1$,
and are bounded in magnitude by $K$ almost surely for some $K $, where
$K \ge10( \E|\xi|^{4} +1)$. Let $H$ be a subspace of dimension $d$ and
$\pi_H$ the orthogonal projection onto $H$. Then
\[
\P\bigl(\bigl|\|\pi_H (X)\| - \sqrt d\bigr| \ge t\bigr) \le10 \exp\biggl(-
\frac{t^{2}}{10K^2}\biggr).
\]
In particular, one has
\[
\| \pi_H(X)\| = \sqrt{d} + O( K \log n )
\]
with overwhelming probability.
\end{lemma}

\begin{pf} See \cite{TVlocal1}, Lemma 43; the proof is a short
application of Talagrand's inequality \cite{Le}.
\end{pf}

\section{Delocalization}

The purpose of this section is to establish Theorem~\ref{delocal} and
Theorem \ref{eigenconc}. The material here is closely analogous to
\cite{TVlocal1}, Sections~5.2, 5.3, as well as that of the original
results in \cite{ESY1,ESY2,ESY3} and can be read independently of the
other sections of the paper. The recent paper \cite{ESYY} also
contains arguments and results closely related to those in this section.

\subsection{\texorpdfstring{Deduction of Theorem \protect\ref{delocal} from Theorem \protect\ref{eigenconc}}{Deduction of Theorem 19 from Theorem 20}}

We begin by showing how Theorem \ref{delocal}\vadjust{\goodbreak} follows from Theorem
\ref{eigenconc}. We shall just establish the claim for the right
singular vectors $u_i$, as the claim for the left singular vectors is
similar. We fix $\eps$ and allow all implied constants to depend on
$\eps$ and $y$. We can also assume that $K^2 \log^{20} n = o(n)$ as
the claim is trivial otherwise.

As $M$ is continuous, we see that the nontrivial singular values are
almost surely simple and positive, so that the singular vectors $u_i$
are well defined up to unit phases. Fix $1 \leq i \leq p$; it suffices
by the union bound and symmetry to show that the event that $\lambda
_i$ falls outside $[a+\eps,b-\eps]$ or that the $n${th} coordinate
$x$ of $u_i$ is $O(K n^{-1/2} \log^{10} n)$ holds with (uniformly)
overwhelming probability.

Applying Corollary \ref{singcoord}, it suffices to show that with
uniformly overwhelming probability, either $\lambda_i \notin[a+\eps
,b-\eps]$, or
%
\begin{equation}\label{jap}
 \qquad \sum_{j=1}^{\min(p,n-1)} \frac{\sigma_j(M_{p,n-1})^2}{(\sigma
_j(M_{p,n-1})^2-\sigma_i(M_{p,n})^2)^{2}} |v_j(M_{p,n-1})^\ast X|^2
\gg\frac{n}{K^2 \log^{20} n},
\end{equation}
where $M =
\left({ M_{p,n-1} \enskip X
}\right)
$. But if $\lambda_i \in[a+\eps,b-\eps]$, then by\footnote{In the
case $p=n$, one would have to replace $M_{p,n-1}$ by its transpose to
return to the regime $p\leq n$.} Theorem \ref{eigenconc}, one can find
(with uniformly overwhelming probability) a set $J \subset\{1,\ldots
,\min(p,\allowbreak n-1)\}$ with $|J| \gg K^2 \log^{20} n$ such that $\lambda
_j(M_{p,n-1}) = \lambda_i(M_{p,n}) +\break O( K^2 \log^{20} n / n )$ for
all $j \in J$; since $\lambda_i = \frac{1}{n} \sigma_i^2$, we
conclude that $\sigma_j(M_{p,n-1})^2 = \sigma_i(M_{p,n})^2 + O( K^2
\log^{20} n )$. In particular, $\sigma_j(M_{p,n-1}) = \Theta(\sqrt
{n})$. By Pythagoras' theorem, the left-hand side of \eqref{jap} is
then bounded from below by
\[
\gg n \frac{\| \pi_H X \|^2}{(K^2 \log^{20} n)^2},
\]
where $H \subset\C^p$ is the span of the $v_j(M_{p,n-1})$ for $j \in
J$. But from Lemma \ref{lemmaprojection} (and the fact that $X$ is
independent of $M_{p,n-1}$), one has
\[
\| \pi_H X \|^2 \gg K^2 \log^{20} n
\]
with uniformly overwhelming probability, and the claim follows.

It thus remains to establish Theorem \ref{eigenconc}.

\subsection{A crude upper bound}

Let the hypotheses be as in Theorem \ref{eigenconc}. We first
establish a crude upper bound, which illustrates the techniques used to
prove Theorem \ref{eigenconc}, and also plays an important direct role
in that proof.

\begin{proposition}[(Eigenvalue upper bound)]\label{eigen-upper} Let the
hypotheses be as in Theorem \ref{delocal}. Then for any interval $I
\subset[a+\eps, b-\eps]$ of length $|I| \geq K \log^{2} n / n$, one
has with overwhelming probability (uniformly in $I$) that
\[
|N_I| \ll n |I|,
\]
where $|I|$ denotes the length of $I$, and $N_I$ was defined in \eqref{NI-def}.
\end{proposition}

To prove this proposition, we suppose for contradiction that
%
\begin{equation}\label{nic}
|N_I| \geq C n |I|
\end{equation}
for some large constant $C$ to be chosen later. We will show that for
$C$ large enough, this leads to a contradiction with overwhelming probability.

We follow the standard approach (see, e.g., \cite{BS}) of controlling
the eigenvalue counting function $N_I$ via the Stieltjes transform
\[
s(z) := \frac{1}{p} \sum_{j=1}^p \frac{1}{\lambda_j(W) - z}.
\]
Fix $I$. If $x$ is the midpoint of $I$, $\eta:= |I|/2$, and $z := x +
\sqrt{-1} \eta$, we see that
\[
\Im s(z) \gg\frac{|N_I|}{\eta p}
\]
[recall that $p = \Theta(n)$] so from \eqref{nic} one has
%
\begin{equation}\label{saz-large}
\Im(s(z)) \gg C.
\end{equation}

Applying Lemma \ref{stielt}, with $W$ replaced by the $p \times p$
matrix $\tilde W := \frac{1}{n} M M^\ast$ (which only has the
nontrivial eigenvalues), we see that
%
\begin{equation}\label{star}
s(z) = \frac{1}{p} \sum_{k=1}^p \frac{1}{\xi_{kk} - z - a_k^\ast
(W_k - zI)^{-1} a_k},
\end{equation}
where $\xi_{kk}$ is the $kk$ entry of $\tilde W$, $W_k$ is the $p-1
\times p-1$ matrix with the $k$th row and column of
$\tilde W$ removed, and $a_k \in\C^{p-1}$ is the $k${th} column of $\tilde W$ with the $k$th entry removed.

Using the crude bound $|\Im\frac{1}{z}| \leq\frac{1}{|\Im(z)|}$
and \eqref{saz-large}, one concludes
\[
\frac{1}{p} \sum_{k=1}^p \frac{1}{|\eta+ \Im a_k^\ast(W_k -
zI)^{-1} a_k|} \gg C.
\]
By the pigeonhole principle, there exists $1 \leq k \leq p$ such that
%
\begin{equation}\label{saz-2}
\frac{1}{|\eta+ \Im a_k^\ast(W_k - zI)^{-1} a_k|} \gg C.
\end{equation}
The fact that $k$ varies will cost us a factor of $p$ in our failure
probability estimates, but this will not be of concern since all of our
claims will hold with overwhelming probability.

Fix $k$. Note that
%
\begin{equation}\label{amk}
a_k = \frac{1}{n} M_k X_k
\end{equation}
and
\[
W_k = \frac{1}{n} M_k M_k^\ast,\vadjust{\goodbreak}
\]
where $X_k \in\C^n$ is the (adjoint of the) $k$th
row of $M$, and $M_k$ is the $p-1 \times n$ matrix formed by removing
that row. Thus, if we let $v_1(M_k),\ldots,v_{p-1}(M_k) \in\C^{p-1}$
and $u_1(M_k),\ldots,u_{p-1}(M_k) \in\C^n$ be coupled orthonormal
systems of left and right singular vectors of $M_k$, and let $\lambda
_j(W_k) = \frac{1}{n} \sigma_j(M_k)^2$ for $1 \leq j \leq p-1$ be the
associated eigenvectors, one has
%
\begin{equation}\label{star2}
a_k^\ast(W_k - zI)^{-1} a_k = \sum_{j=1}^{p-1} \frac{|a_k^\ast
v_j(M_k)|^2}{\lambda_j(W_k) - z}.
\end{equation}
and thus
\[
\Im a_k^\ast(W_k - zI)^{-1} a_k \geq\eta\sum_{j=1}^{p-1} \frac
{|a_k^\ast v_j(M_k)|^2}{\eta^2 + |\lambda_j(W_k) - x|^2}.
\]
We conclude that
\[
\sum_{j=1}^{p-1} \frac{|a_k^\ast v_j(M_k)|^2}{\eta^2 + |\lambda
_j(W_k) - x|^2} \ll\frac{1}{C\eta}.
\]

The expression $a_k^\ast v_j(M_k)$ can be rewritten much more favorably
using~\eqref{amk}~as
%
\begin{equation}\label{star3}
a_k^\ast v_j(M_k) = \frac{\sigma_j(M_k)}{n} X_k^\ast u_j(M_k).
\end{equation}
The advantage of this latter formulation is that the random variables
$X_k$ and $u_j(M_k)$ are independent (for fixed $k$).

Next, note that from \eqref{nic} and the Cauchy interlacing law (Lemma
\ref{cauchy}) one can find an interval $J \subset\{1,\ldots,p-1\}$
of length
%
\begin{equation}\label{jeta}
|J| \gg C \eta n
\end{equation}
such that $\lambda_j(W_k) \in I$. We conclude that
\[
\sum_{j \in J} \frac{\sigma_j(M_k)^2}{n^2} |X_k^\ast u_j(M_k)|^2 \ll
\frac{\eta}{C}.
\]
Since $\lambda_j(W_k) \in I$, one has $\sigma_j(M_k) = \Theta(\sqrt
{n})$, and thus
\[
\sum_{j \in J} |X_k^\ast u_j(M_k)|^2 \ll\frac{\eta n}{C}.
\]
The left-hand side can be rewritten using Pythagoras' theorem as $\|
\pi_H X_k \|^2$, where $H$ is the span of the eigenvectors $u_j(M_k)$
for $j \in J$. But from Lemma~\ref{lemmaprojection} and \eqref{jeta}, we see that this quantity is $\gg\eta n$ with overwhelming
probability, giving the desired contradiction with overwhelming
probability (even after taking the union bound in~$k$). This concludes
the proof of Proposition \ref{eigen-upper}.

\subsection{Reduction to a Stieltjes transform bound}

We now begin the proof of Theorem \ref{eigenconc} in earnest. We
continue to allow all implied constants to depend on $\eps$
and~$y$.\vadjust{\goodbreak}

It suffices by a limiting argument (using Lemma \ref{weyl}) to
establish the claim under the assumption that the distribution of $M$
is continuous; our arguments will not use any quantitative estimates on
this continuity.

The strategy is to compare $s$ with the Marchenko--Pastur Stieltjes transform
\[
s_{\mathrm{MP},y}(z) := \int_\R\rho_{\mathrm{MP},y}(x) \frac{1}{x-z}\,dx.
\]
A routine application of \eqref{rhomp} and the Cauchy integral formula
yields the explicit formula
%
\begin{equation}\label{lemon}
s_{\mathrm{MP},y}(z) = - \frac{y+z-1- \sqrt{(y+z-1)^2 - 4yz}}{2yz},
\end{equation}
where we use the branch of $\sqrt{(y+z-1)^2 - 4yz}$ with cut at
$[a,b]$ that is asymptotic to $y-z+1$ as $z \to\infty$. To put it
another way, for $z$ in the upper half-plane, $s_{\mathrm{MP},y}(z)$ is the
unique solution to the equation
%
\begin{equation}\label{smpy}
s_{\mathrm{MP},y} = - \frac{1}{y+z-1 + yz s_{\mathrm{MP},y}(z)}
\end{equation}
with $\Im s_{\mathrm{MP},y}(z) > 0$. (Details of these computations can also be
found in~\cite{BS}.)

We have the following standard relation between convergence of
Stieltjes transform and convergence of the counting function.

\begin{lemma}[(Stieltjes transform controls counting function)]\label{lemmaS-transform}
Let $1/10 \geq\eta\geq1/n$, and $L, \eps, \delta> 0$. Suppose that
one has the bound
%
\begin{equation}\label{soda}
|s_{\mathrm{MP},y}(z) -s(z) | \le\delta
\end{equation}
with overwhelming probability for each $z$ with $|\Re(z)| \leq L$ and
$\Im(z) \geq\eta$, with the implied constants in the definition of
overwhelming probability uniform in $z$. Then for any interval $I$ in
$[a+\eps,b-\eps]$ with $|I| \geq\max( 2\eta, \frac{\eta}{\delta
} \log\frac{1}{\delta} )$, one has
\[
\biggl|N_I - n \int_I \rho_{\mathrm{MP},y}(x)\,dx\biggr| \ll\delta n |I|
\]
with overwhelming probability.
\end{lemma}

\begin{pf} This follows from \cite{TVlocal1}, Lemma 64; strictly
speaking, that lemma was phrased for the semi-circular distribution
rather than the Marchenko--Pastur distribution, but an inspection of the
proof shows the proof can be modified without difficulty. See also
\cite{GT00} and \cite{ESY1}, Corollary 4.2, for closely related lemmas.
\end{pf}

In view of this lemma, we see that to show Theorem \ref{eigenconc}, it
suffices to show that for each complex number $z$ in the region
\[
\Omega:= \biggl \{ z \in\C\dvtx  a+\eps/2 \leq\Re(z) \leq b-\eps/2;
 \Im(z) \geq \eta:= \frac{K^2 \log^{19} n}{n}  \biggr\},\vadjust{\goodbreak}
\]
one has
%
\begin{equation}\label{sang}
s(z) - s_{\mathrm{MP},y}(z) = o(1)
\end{equation}
with (uniformly) overwhelming probability.

For this, we return to the formula \eqref{star}. Inserting the
identities \eqref{star2}, \eqref{star3} into this formula, one obtains
%
\begin{equation}\label{start}
s(z) = \frac{1}{p} \sum_{k=1}^p \frac{1}{\xi_{kk} - z - Y_k},
\end{equation}
where $Y_k = Y_k(z)$ is the quantity
\[
Y_k := \sum_{j=1}^{p-1} \frac{\lambda_j(M_k)}{n} \frac{ |X_k^\ast
u_j(M_k)|^2}{\lambda_j(W_k) - z}.
\]
Suppose we condition $M_k$ (and thus $W_k$) to be fixed; the entries of
$X_k$ remain independent with mean zero and variance $1$, and thus
(since the $u_j$ are unit vectors)
\begin{eqnarray*}
\E(Y_k|M_k)
&=& \sum_{j=1}^{p-1} \frac{\lambda_j(M_k)}{n} \frac{
1}{\lambda_j(W_k) - z} \\
&=& \frac{p-1}{n} \bigl(1 + z s_k(z) \bigr),
\end{eqnarray*}
where
\[
s_k(z) := \frac{1}{p-1} \sum_{j=1}^{p-1} \frac{1}{\lambda_j(W_k) -
z}
\]
is the Stieltjes transform of $W_k$.

From the Cauchy interlacing law (Lemma \ref{cauchy}), we see that the
difference
\[
s(z) - \frac{p-1}{p} s_k(z) = \frac{1}{p}  \Biggl( \sum_{j=1}^p \frac
{1}{\lambda_j(W) - z} - \sum_{j=1}^{p-1} \frac{1}{\lambda_j(W_k) -
z}  \Biggr)
\]
is bounded in magnitude by $O(\frac{1}{p})$ times the total variation
of the function $\lambda\mapsto\frac{1}{\lambda-z}$ on $[0,+\infty
)$, which is $O( \frac{1}{\eta} )$. Thus,
\[
\frac{p-1}{p} s_k(z) = s(z) + O \biggl( \frac{1}{p \eta}  \biggr)
\]
and thus
%
\begin{eqnarray}\label{ykmk}
\E( Y_k | M_k ) &=& \frac{p-1}{n} + \frac{p}{n} z s(z) + O \biggl(
\frac{1}{n \eta}  \biggr) \nonumber
\\[-8pt]
\\[-8pt]
&=& y + o(1) + \bigl(y+o(1)\bigr) z s(z)
\nonumber
\end{eqnarray}
since $p/n= y+o(1)$ and $1/\eta= o(n)$.\vadjust{\goodbreak}

We will shortly show a similar bound for $Y_k$ itself.\vspace*{-3pt}

\begin{lemma}[(Concentration of $Y_k$)]\label{conc} Let $z \in\Omega$.
For each $1 \leq k \leq p$, one has $Y_k = y + o(1) + (y+o(1)) z s(z)$
with overwhelming probability (uniformly in $k$ and $I$).\vspace*{-3pt}
\end{lemma}

Meanwhile, we have
\[
\xi_{kk} = \frac{1}{n} \| X_k \|^2
\]
and hence by Lemma \ref{lemmaprojection}, $\xi_{kk} = 1+o(1)$ with
overwhelming probability (again uniformly in $k$ and $I$). Inserting
these bounds into \eqref{start}, one obtains
\[
s(z) = \frac{1}{p} \sum_{k=1}^p \frac{1}{1 - z - (y+o(1)) - (y+o(1))
z s(z)}
\]
with overwhelming probability; thus $s(z)$ ``almost solves'' \eqref{smpy} in some sense. From the quadratic formula, the two solutions of
\eqref{smpy} are $s_{\mathrm{MP},y}(z)$ and $-\frac{y+z-1}{yz} - s_{\mathrm{MP},y}(z)$.
One concludes that for each fixed $z \in\Omega$, it occurs with
overwhelming probability that one has either
%
\begin{equation}\label{saz}
s(z) = s_{\mathrm{MP},y}(z) + o(1)
\end{equation}
or
%
\begin{equation}\label{saz2}
s(z) = - \frac{y+z-1}{yz} + o(1)
\end{equation}
or
%
\begin{equation}\label{saz3}
s(z) = - \frac{y+z-1}{yz} - s_{\mathrm{MP},y}(z) + o(1)
\end{equation}
(with the convention that $\frac{y+z-1}{yz} = 1$ when $y=1$). By using
a $n^{-100}$-net of possible $z$'s in $\Omega$ and using the union
bound [and the fact that $s(z)$ has a Lipschitz constant of at most
$O(n^{10})$ in $\Omega$] we may assume (with overwhelming probability)
that the above trichotomy holds for \textit{all} $z \in\Omega$. In
other words, if $\delta> 0$ is a small number (which may depend on
$a,b,\eps$) and $n$ is sufficiently large depending on $\delta$, we
may cover
\[
\Omega\subset\Omega_1 \cup\Omega_2 \cup\Omega_3,
\]
where
\begin{eqnarray*}
\Omega_1 &:=& \{ z \in\Omega\dvtx  |s(z) - s_{\mathrm{MP},y}(z)| \leq\delta\},\\[-2pt]
\Omega_2 &:=& \biggl\{ z \in\Omega\dvtx  \biggl|s(z) + \frac{y+z-1}{yz}\biggr| \leq\delta
\biggr\},\\[-2pt]
\Omega_3 &:=& \biggl\{ z \in\Omega\dvtx  \biggl|s(z) + \frac{y+z-1}{yz} +
s_{\mathrm{MP},y}(z)\biggr| \leq\delta\biggr\}.
\end{eqnarray*}

When $\Im(z) = n^{10}$, then $s(z), s_{\mathrm{MP},y}(z)$ are both $o(1)$, and
so (for $n$ sufficiently large) we see that $z \in\Omega_1$ in this
case. In particular, $\Omega_1$ is empty.\vadjust{\goodbreak} On the other hand, $\Omega
_1, \Omega_2, \Omega_3$ are closed subsets of $\Omega$.
From \eqref{smpy}, one has
\[
s_{\mathrm{MP},y}(z)  \biggl( \frac{y+z-1}{yz} + s_{\mathrm{MP},y}(z)  \biggr) = -
\frac{1}{yz},
\]
which implies that the separation between $s_{\mathrm{MP},y}(z)$ from $-\frac
{y+z-1}{yz}$ is bounded from below, which implies that $\Omega_1$ and
$\Omega_2$ are disjoint (for $\delta$ small enough). Similarly, from
\eqref{lemon}, we see that
\[
\frac{y+z-1}{yz} + 2 s_{\mathrm{MP},y}(z) = \frac{\sqrt{(y+z-1)^2 - 4yz}}{yz};
\]
since $(y+z-1)^2-4yz$ has zeroes only when $z=a,b$, and $z$ is bounded
away from these singularities, we see also that $\Omega_1$ and $\Omega
_3$ are also disjoint.

The sets $\Omega_1$, $\Omega_2 \cup\Omega_3$ are thus disjoint
closed subsets of $\Omega$. As $\Omega$ is connected and $\Omega_1$
is nonempty, we conclude that $\Omega_1 = \Omega$ (whenever $n$ is
sufficiently large depending on $\delta$). Letting $\delta\to0$, we
conclude that \eqref{saz} holds unniformly for $z \in\Omega$ with
overwhelming probability, which gives \eqref{sang} and thus Theorem
\ref{eigenconc}.

\section{\texorpdfstring{Proof of Theorem \protect\ref{four-main} and Theorem \protect\ref{gap1}}{Proof of Theorem 17 and Theorem 18}}\label{fmt}

We first prove Theorem \ref{four-main}. The arguments follow those in
\cite{TVlocal1}.

We begin by observing from Markov's inequality and the union bound that
one has $|\zeta_{ij}|, |\zeta'_{ij}| \leq n^{10/C_0}$ (say) for all
$i,j$ with probability $O( n^{-8})$. Thus, by truncation (and adjusting
the moments appropriately, using Lemma \ref{weyl} to absorb the
error), one may assume without loss of generality that
%
\begin{equation}\label{zetas}
|\zeta_{ij}|, |\zeta'_{ij}| \leq n^{10/C_0}
\end{equation}
almost surely for all $i,j$.
Next, by a further approximation argument we may assume that the
distribution of $M, M'$ is continuous. This is a purely qualitative
assumption, to ensure that the singular values are almost surely
simple; our bounds will not depend on any quantitative measure on the
continuity, and so the general case then follows by a limiting argument
using Lemma \ref{weyl}.

The key technical step is the following theorem, whose proof is delayed
to the next section.

\begin{theorem}[(Truncated Four Moment theorem)]\label{four-main-truncate}
For sufficiently small $c_0>0$ and sufficiently large $C_0>0$, the
following holds for every $0 < \eps< 1$ and $k \geq1$.
Let $M = (\zeta_{ij})_{1 \leq i \leq p, 1 \leq j \leq n}$ and $M' =
(\zeta'_{ij})_{1 \leq i \leq p, 1 \leq j \leq n}$ be matrix ensembles
obeying condition \textup{C1} for some $C_0$, as well as \eqref{zetas}.
Assume that $p/n \to y$ for some $0 < y \leq1$, and that $\zeta_{ij}$
and $\zeta'_{ij}$ match to order $4$.

Let $G\dvtx  \R^k \times\R_+^k \to\R$ be a smooth function obeying the
derivative bounds
%
\begin{equation}\label{G-deriv-2}
|\nabla^j G(x_1,\ldots,x_k,q_1,\ldots,q_k)| \leq n^{c_0}
\end{equation}
for all $0 \leq j \leq5$ and $x_1,\ldots,x_k \in\R$, $q_1,\ldots
,q_k \in\R$, and such that $G$ is supported on the region $q_1,\ldots
,q_k \leq n^{c_0}$, and the gradient $\nabla$ is in all $2k$ variables.

Then for any $\eps p \le i_1 < i_2 <\cdots< i_k \le(1-\eps) p$, and
for $n$ sufficiently large depending on $\eps, k, c_0$ we have
%
\begin{eqnarray} \label{eqnapproximation-2}
&& \quad  \bigl|\E\bigl( G\bigl(\sqrt{n} \sigma_{i_1}(M), \ldots , \sqrt{n} \sigma
_{i_k}(M), Q_{i_1}(M), \ldots,Q_{i_k}(M)\bigr)\bigr) \nonumber
\\[-8pt]
\\[-8pt]
&& \quad \hphantom{\bigl|}{} -
\E\bigl( G\bigl(\sqrt{n} \sigma_{i_1}(M'), \ldots , \sqrt{n} \sigma
_{i_k}(M'), Q_{i_1}(M'), \ldots,Q_{i_k}(M')\bigr)\bigr)\bigr| \le n^{-c_0}.
\nonumber
\end{eqnarray}
If $\zeta_{ij}, \zeta'_{ij}$ match to order $3$, then the conclusion
still holds as long as one strengthens \eqref{G-deriv-2} to
%
\begin{equation}\label{rock}
|\nabla^j G(x_1,\ldots,x_k,q_1,\ldots,q_k)| \leq n^{-jc_1}
\end{equation}
for some $c_1>0$, if $c_0$ is sufficiently small depending on $c_1$.
\end{theorem}

Informally, Theorem \ref{four-main-truncate} is a truncated version of
Theorem \ref{four-main} in which one has smoothly restricted attention
to the event where eigenvalue gaps are not unexpectedly small.

Given a $p \times n$ matrix $M$ we form the augmented matrix $\M$
defined in \eqref{augment}, whose eigenvalues are $\pm\sigma
_1(M),\ldots,\pm\sigma_p(M)$, together with the eigenvalue~$0$ with
multiplicity $n-p$ (if $p<n$). For each $1 \leq i \leq p$, we introduce
(in analogy with the arguments in \cite{TVlocal1}) the quantities
\begin{eqnarray*}
&&Q_i(\M)\\
 && \qquad := \sum_{\lambda\neq\sigma_i(M)} \frac{1}{|\sqrt{n}(
\lambda- \sigma_i(M) )|^2} \\
&& \qquad \hspace*{3pt}= \frac{1}{n} \Biggl( \sum_{1 \leq j \leq p:  j \neq i} \frac
{1}{|\sigma_j(M) -\sigma_i(M)|^2} + \frac{n-p}{\sigma_i(M)^2} +
\sum_{j=1}^p \frac{1}{|\sigma_j(M) + \sigma_i(M)|^2}  \Biggr).
\end{eqnarray*}
(The factor of $\frac{1}{n}$ in $Q_i(\M)$ is present to align the
notation here with that in \cite{TVlocal1}, in which one dilated the
matrix by $\sqrt{n}$.) We set $Q_i(\M)=\infty$ if the singular value
$\sigma_i$ is repeated, but this event occurs with probability zero
since we are assuming $M$ to be continuously distributed. One should
view~$Q_i(\M)$ as measuring the extent to which eigenvalue (or
singular value) gaps near~$\sigma_i(M)$ are unexpectedly small.

The gap property on $M$ ensures an upper bound on $Q_i(\M)$.

\begin{lemma}\label{highprob} If $M$ satisfies the gap property, then
for any $c_0 > 0$ (independent of $n$), and any $\eps p \leq i \leq
(1-\eps) p$, one has $Q_i(\M) \leq n^{c_0}$ with high probability.
\end{lemma}

\begin{pf} Observe the upper bound
%
\begin{equation}\label{qi}
Q_i(\M) \leq\frac{2}{n} \sum_{1 \leq j \leq p:  j \neq i} \frac
{1}{|\sigma_j(M) -\sigma_i(M)|^2} +
\frac{n-p+1}{n\sigma_i(M)^2}.\vadjust{\goodbreak}
\end{equation}

From Corollary \ref{bulk}, we see that with overwhelming probability,
$\sigma_i(M)^2 / n$ is bounded away from zero, and so $\frac
{n-p+1}{n\sigma_i(M)^2} = O(1/n)$. To bound the other term in \eqref{qi}, one
repeats the proof of \cite{TVlocal1}, Lemma 49.\vspace*{-3pt}
\end{pf}


By applying a truncation argument exactly as in \cite{TVlocal1}, Section 3.3, one can now remove the hypothesis in Theorem \ref{four-main-truncate} that $G$ is supported in the region $q_1,\ldots
,q_k \leq n^{c_0}$. In particular, one can now handle the case when $G$
is independent of $q_1,\ldots,q_k$; and Theorem \ref{four-main}
follows after making the change of variables $\lambda= \frac{1}{n}
\sigma^2$ and using the chain rule (and Corollary \ref{bulk}).

Next, we prove Theorem \ref{gap1}, assuming both Theorems \ref{four-main-truncate} and \ref{gap}.
The main observation here is the following lemma.\vspace*{-3pt}


\begin{lemma}[(Matching lemma)] \label{lemmamatching} Let $\zeta$ be a
complex random variable with mean zero, unit variance, and third moment
bounded by some constant
$a$. Then there exists a complex random variable $\tilde\zeta$ with
support bounded by the ball of radius $O_a(1)$ centered at the origin
(and in particular, obeying the exponential decay hypothesis uniformly
in $\zeta$ for fixed $a$) which matches $\zeta$ to third order.\vspace*{-3pt}
\end{lemma}

\begin{pf} In order for $\tilde\zeta$ to match $\zeta$ to third
order, it suffices that $\tilde\zeta$ have mean zero, variance $1$,
and that $\E\tilde\zeta^3 = \E\zeta^3$ and $\E\tilde\zeta^2
\overline{\tilde\zeta} = \E\zeta^3 \overline{\zeta}$.

Accordingly, let $\Omega\subset\C^2$ be the set of pairs $(\E\tilde
\zeta^3, \E\tilde\zeta^2 \overline{\tilde\zeta})$ where $\tilde
\zeta$ ranges over complex random variables with mean zero, variance
one, and compact support. Clearly $\Omega$ is convex. It is also
invariant under the symmetry $(z,w) \mapsto(e^{3i\theta} z,
e^{i\theta} w)$ for any phase $\theta$. Thus, if $(z,w) \in\Omega$,
then $(-z,e^{i \pi/3} w) \in\Omega$, and hence by convexity
$(0,\frac{\sqrt{3}}{2} e^{i \pi/6} w) \in\Omega$, and hence by
convexity and rotation invariance $(0,w') \in\Omega$ whenever $|w'|
\leq\frac{\sqrt{3}}{2} w$. Since $(z,w)$ and $(0,-\frac{\sqrt
{3}}{2} w)$ both lie in $\Omega$, by convexity $(cz,0)$ lies in it
also for some absolute constant $c>0$, and so again by convexity and
rotation invariance $(z',0) \in\Omega$ whenever $|z'| \leq cz$. One
last application of convexity then gives $(z'/2,w'/2) \in\Omega$
whenever $|z'| \leq cz$ and $|w'| \leq\frac{\sqrt{3}}{2} w$.

It is easy to construct complex random variables with mean zero,
variance one, compact support, and arbitrarily large third moment.
Since the third moment is comparable to $|z|+|w|$, we thus conclude
that $\Omega$ contains all of~$\C^2$, that is, every complex random
variable with finite third moment with mean zero and unit variance can
be matched to third order by a variable of compact support. An
inspection of the argument shows that if the third moment is bounded by
$a$ then the support can also be bounded by $O_a(1)$.\vspace*{-3pt}~%
\end{pf}

Now consider a random matrix $M$ as in Theorem \ref{gap1} with atom
variables~$\zeta_{ij}$.
By the above lemma, for each $i,j$, we can find $ \zeta'_{ij}$ which
satisfies the exponential decay hypothesis and match
$\zeta_{ij}$ to third order. Let $\eta(q)$ be a~smooth cutoff to the
region $q \leq n^{c_0}$ for some $c_0>0$ independent of $n$, and let\vadjust{\goodbreak}
$\eps p \leq i \leq(1-\eps) p$.
By Theorem \ref{gap}, the matrix $M'$ formed by the $\zeta'_{ij}$
satisfies the gap property. By Lemma~\ref{highprob},
\[
\E \eta( Q_i(M') ) = 1 - O( n^{-c_1} )
\]
for some $c_1 > 0$ independent of $n$, so by Theorem \ref{four-main-truncate} one has
\[
\E \eta( Q_i(M) ) = 1 - O( n^{-c_2} )
\]
for some $c_2 > 0$ independent of $n$. We conclude that $M$ also obeys
the gap property.

The next two sections are devoted to the proofs of Theorem \ref{four-main-truncate} and Theorem~\ref{gap}, respectively.

\begin{remark}\label{remove} The above trick to remove the exponential
decay hypothesis for Theorem \ref{gap} also works to remove the same
hypothesis in \cite{TVlocal1}, Theorem 19. The point is that in the
analogue of Theorem \ref{four-main-truncate} in that paper (implicit
in \cite{TVlocal1}, Section~3.3), the exponential decay hypothesis is
not used anywhere in the argument; only a uniformly bounded $C_0$
moment for $C_0$ large enough is required, as is the case here. Because
of this, one can replace all the exponential decay hypotheses in the
results of \cite{TVlocal1,TVlocal2} by a hypothesis of bounded $C_0$
moment; we omit the details.\vspace*{-3pt}
\end{remark}

\section{\texorpdfstring{The proof of Theorem \protect\ref{four-main-truncate}}{The proof of Theorem 32}} \label{sectionfmtr}

It remains to prove Theorem \ref{four-main-truncate}. By telescoping
series, it suffices to establish a bound
%
\begin{eqnarray} \label{eqnapproximation-3}
&&\bigl|\E\bigl( G\bigl(\sqrt{n} \sigma_{i_1}(M), \ldots , \sqrt{n} \sigma
_{i_k}(M), Q_{i_1}(M), \ldots,Q_{i_k}(M)\bigr)\bigr) \nonumber\\[-2pt]
&&\hphantom{\bigl|} {}-
\E\bigl( G\bigl(\sqrt{n} \sigma_{i_1}(M'), \ldots , \sqrt{n} \sigma
_{i_k}(M'), Q_{i_1}(M'), \ldots,Q_{i_k}(M')\bigr)\bigr)\bigr| \\[-2pt]
&& \qquad   \le n^{-2-c_0}\nonumber
\end{eqnarray}
under the assumption that the coefficients $\zeta_{ij}$, $\zeta
'_{ij}$ of $M$ and $M'$ are identical except in one entry, say the $qr$
entry for some $1 \leq q \leq p$ and $1 \leq r \leq n$, since the claim
then follows by interchanging each of the $pn = O(n^2)$ entries of $M$
into $M'$ separately.

Write $M(z)$ for the matrix $M$ (or $M'$) with the $qr$ entry replaced
by $z$. We apply the following proposition, which follows from a
lengthy argument in \cite{TVlocal1}:\vspace*{-3pt}

\begin{proposition}[(Replacement given a good configuration)]\label{swap}
Let the notation and assumptions be as in Theorem \ref{four-main-truncate}.
There exists a positive constant~$C_1$
(independent of $k$) such that the following holds. Let $\eps_1 > 0$.
We condition (i.e., freeze) all the entries of $M(z)$ to be constant,
except for the $qr$ entry, which is $z$. We assume that for every $1
\leq j \leq k$ and every $|z| \leq n^{1/2+\eps_1}$ whose real and
imaginary parts are multiples of $n^{-C_1}$, we have
\begin{itemize}
\item(Singular value separation) For any $1 \leq i \leq n$ with
$|i-i_j| \geq n^{\eps_1}$, we have
%
\begin{equation}\label{noon}
\bigl|\sqrt{n} \sigma_i(M(z)) - \sqrt{n} \sigma_{i_j}(M(z))\bigr| \geq
n^{-\eps_1} |i-i_j|.
\end{equation}
Also, we assume
%
\begin{equation}\label{noon2}
\sqrt{n} \sigma_{i_j}(A(z)) \geq n^{-\eps_1} n.\vadjust{\goodbreak}
\end{equation}
\item(Delocalization at $i_j$) If $u_{i_j}(M(z)) \in\C^n$,
$v_{i_j}(M(z)) \in\C^p$ are unit right and left singular vectors of
$M(z)$, then
%
\begin{equation}\label{pz1}
|e_q^* v_{i_j}(M(z))|, |e_r^* u_{i_j}(M(z))| \leq n^{-1/2+\eps_1}.
\end{equation}
\item For every $\alpha\geq0$
%
\begin{equation}\label{pz2}
\| P_{i_j,\alpha}(M(z)) e_q \|, \| P'_{i_j,\alpha}(M(z)) e_r \| \leq
2^{\alpha/2} n^{-1/2+\eps_1},
\end{equation}
whenever $P_{i_j,\alpha}$ (resp., $P'_{i_j,\alpha}$) is the orthogonal
projection to the span of right singular vectors $u_i(M(z))$ [resp.,
left singular vectors $v_i(M(z))$] corresponding to singular values
$\sigma_i(A(z))$ with $2^\alpha\leq|i-i_j| < 2^{\alpha+1}$.
\end{itemize}
We say that $M(0), e_q, e_r$ are a \textit{good configuration} for
$i_1,\ldots,i_k$ if the above properties hold. Assuming this good
configuration, then we have \eqref{eqnapproximation-3} if $\zeta
_{ij}$ and $\zeta'_{ij}$ match to order $4$, or if they match to order
$3$ and \eqref{rock} holds.\vspace*{-3pt}
\end{proposition}

\begin{pf}
$\!\!\!$This follows by applying \cite{TVlocal1}, Proposition 46,
to the \mbox{$p+n \times p+n$} Hermitian matrix $A(z) := \sqrt{n} \M(z)$,
where $\M(z)$ is the augmented matrix of~$M(z)$, defined in \eqref{augment}. Note that the eigenvalues of $A(z)$ are $\pm\sqrt{n}
\sigma_1(M(z)),\break \ldots, \pm\sqrt{n} \sigma_p(M(z))$ and $0$, and
that the eigenvalues are given (up to unit phases) by $
{ v_j(M(z)) \choose \pm u_j(M(z))
}
$. Note also that the analogue of \eqref{pz2} in \cite{TVlocal1},
Proposition 46, is trivially true if~$2^\alpha$ is comparable to $n$, so
one can restrict attention to the regime $2^\alpha= o(n)$.\vspace*{-3pt}
\end{pf}

In view of the above proposition, we see that to conclude the proof of
Theorem~\ref{four-main-truncate} (and thus Theorem \ref{four-main})
it suffices to show that for any $\eps_1 > 0$, that~$M(0)$, $e_q$, $e_r$
are a good configuration for $i_1,\ldots,i_k$ with overwhelming
probability, if $C_0$ is sufficiently large depending on $\eps_1$ (cf.
\cite{TVlocal1}, Proposition~48).\looseness=-1

Our main tools for this are Theorem \ref{delocal} and Theorem \ref{eigenconc}. Actually,
we need a slight variant.\vspace*{-3pt}

\begin{proposition}\label{sdb2} The conclusions of Theorem \ref{delocal} and Theorem \ref{eigenconc} continue to hold if one replaces
the $qr$ entry of $M$ by a deterministic number $z = O( n^{1/2 +
O(1/C_0)} )$.\vspace*{-3pt}
\end{proposition}

This is proven exactly as in \cite{TVlocal1}, Corollary 63, and is omitted.

We return to the task of establishing a good configuration with
overwhelming probability.
By the union bound, we may fix $1 \leq j \leq k$, and also fix the $|z|
\leq n^{1/2+\eps_1}$ whose real and imaginary parts are multiples of
$n^{-C_1}$. By the union bound again and Proposition \ref{sdb2}, the
eigenvalue separation condition \eqref{noon} holds with overwhelming
probability for every $1 \leq i \leq n$ with $|i-j| \geq n^{\eps_1}$
(if $C_0$ is sufficiently large), as does \eqref{pz1}. A similar
argument using Pythagoras' theorem and Corollary \ref{bulk} gives
\eqref{pz2} with overwhelming probability [noting as before that we
may restrict attention to the regime $2^\alpha= o(n)$]. Corollary \ref{bulk} also gives \eqref{noon2} with overwhelming probability. This
gives the claim, and Theorem \ref{four-main} follows.\vadjust{\goodbreak}\vspace*{-3pt}

\section{\texorpdfstring{Proof of Theorem \protect\ref{gap}}{Proof of Theorem 16}} \label{sectiongap}

We now prove Theorem \ref{gap}, closely following the analogous
arguments in \cite{TVlocal1}. Using the exponential decay condition,
we may truncate the $\zeta_{ij}$ (and renormalise moments, using Lemma
\ref{weyl}) to assume that\looseness=-1
%
\begin{equation}\label{supi}
|\zeta_{ij}| \leq\log^{O(1)} n
\end{equation}\looseness=0
almost surely. By a limiting argument, we may assume that $M$ has a
continuous distribution, so that the singular values are almost surely simple.

We write $i_0$ instead of $i$, $p_0$ instead of $p$, and write $N_0 :=
p_0+n$. As in~\cite{TVlocal1}, the strategy is to propagate a narrow
gap for $M = M_{p_0,n}$ backwards in the~$p$ variable, until one can
use Theorem \ref{eigenconc} to show that the gap occurs with small probability.

More precisely, for any $1 \leq i-l < i \leq p \leq p_0$, we let
$M_{p,n}$ be the $p \times n$ matrix formed using the first $p$ rows of
$M_{p_0,n}$, and we define (following \cite{TVlocal1}) the \textit
{regularized gap}
%
\begin{equation}\label{giln}
g_{i,l,p} := \inf_{1 \leq i_- \leq i-l < i \leq i_+ \leq p} \frac
{\sqrt{N_0} \sigma_{i_+}(M_{p,n})-\sqrt{N_0} \sigma
_{i_-}(M_{p,n})}{\min( i_+-i_-, \log^{C_1} N_0 )^{\log^{0.9} N_0}},
\end{equation}
where $C_1 > 1$ is a large constant to be chosen later. It will suffice
to show that
%
\begin{equation}\label{goin}
g_{i_0,1,p_0} \leq n^{-c_0}.
\end{equation}

The main tool for this is
the following lemma.\vspace*{-3pt}
\begin{lemma}[(Backwards propagation of gap)]\label{backprop} Suppose
that $p_0/2 \leq p < p_0$ and $l \leq\eps p/10$ is such that
%
\begin{equation}\label{gdel}
g_{i_0,l,p+1} \leq\delta
\end{equation}
for some $0 < \delta\leq1$ (which can depend on $n$), and that
%
\begin{equation}\label{gilp}
g_{i_0,l+1,p} \geq2^m g_{i_0,l,p+1}
\end{equation}
for some $m \geq0$ with
%
\begin{equation}\label{mcivil}
2^m \leq\delta^{-1/2}.
\end{equation}
Let $X_{p+1}$ be the ($p+1$)th row of $M_{p_0,n}$, and
let $u_1(M_{p,n}),\ldots,u_p(M_{p,n})$ be an orthonormal system of
right singular vectors of $M_{p,n}$ associated to $\sigma
_1(M_{p,n}),\ldots, \sigma_p(M_{p,n})$.
Then one of the following statements hold:
\begin{longlist}[(iii)]
\item[(i)] (Macroscopic spectral concentration) There exists \mbox{$1 \leq
i_- < i_+ \leq p+1$} with $i_+-i_- \geq\log^{C_1/2} n$ such that
$|\sqrt{n} \sigma_{i_+}(M_{p+1,n}) - \sqrt{n} \sigma
_{i_-}(M_{p+1,n})| \leq\break\delta^{1/4} \exp( \log^{0.95} n )
(i_+-i_-)$.
\item[(ii)] (Small inner products) There exists $\eps p/2 \leq i_-
\leq i_0-l < i_0 \leq i_+ \leq(1-\eps/2) p$ with $i_+-i_- \leq\log
^{C_1/2} n$ such that
%
\begin{equation}\label{smallin}
\sum_{i_- \leq j < i_+} |X_{p+1}^* u_j(M_{p,n})|^2 \leq\frac
{i_+-i_-}{2^{m/2} \log^{0.01} n}.\vadjust{\goodbreak}
\end{equation}
\item[(iii)] (Large singular value) For some $1 \leq i \leq p+1$, one has
\[
|\sigma_i(M_{p+1,n})| \geq\frac{\sqrt{n} \exp( -\log^{0.95} n
)}{\delta^{1/2}}.
\]
\item[(iv)]  (Large inner product in bulk) There exists $\eps p/10 \leq
i \leq(1-\eps/10) p$ such that
\[
|X_{p+1}^* u_i(M_{p,n})|^2 \geq\frac{\exp( - \log^{0.96} n
)}{\delta^{1/2}}.
\]
\item[(v)]  (Large row) We have
\[
\|X_{p+1}\|^2 \geq\frac{n \exp( - \log^{0.96} n )}{\delta^{1/2}}.
\]
\item[(vi)]  (Large inner product near $i_0$) There exists $\eps p/10
\leq i \leq(1-\eps/10) p$ with $|i-i_0| \leq\log^{C_1} n$ such that
\[
|X_{p+1}^* u_i(M_{p,n})|^2 \geq2^{m/2} n \log^{0.8} n.\vspace*{-3pt}
\]
\end{longlist}
\end{lemma}

\begin{pf} This follows by applying\footnote{Strictly speaking,
there are some harmless adjustments by constant factors that need to be
made to this lemma, ultimately coming from the fact that $n, p, n+p$
are only comparable up to constants, rather than equal, but these
adjustments make only a negligible change to the proof of that lemma.}
\cite{TVlocal1}, Lemma 51, to the $p+n+1 \times p+n+1$ Hermitian matrix
\[
A_{p+n+1} := \sqrt{n}
\pmatrix{\displaystyle 0 & M_{p+1,n}^\ast\cr\displaystyle  M_{p+1,n} & 0
}
,
\]
which after removing the bottom row and rightmost column (which is
$X_{p+1}$, plus $p+1$ zeroes) yields the $p+n \times p+n$ Hermitian matrix
\[
A_{p+n} := \sqrt{n}
\pmatrix{\displaystyle 0 & M_{p,n}^\ast\cr\displaystyle  M_{p,n} & 0
}
\]
which has eigenvalues $\pm\sqrt{n} \sigma_1(M_{p,n}), \ldots, \pm
\sqrt{n} \sigma_p(M_{p,n})$ and $0$, and an orthonormal eigenbasis
that includes the vectors $
{ u_j(M_{p,n}) \choose v_j(M_{p,n})
}
$ for $1 \leq j \leq p$. (The ``large coefficient'' event in \cite{TVlocal1},
Lemma 51(iii), cannot occur here, as $A_{p+n+1}$ has zero diagonal.)\vspace*{-3pt}
\end{pf}

By repeating the arguments in \cite{TVlocal1}, Section 3.5, almost
verbatim, it then suffices to show the following proposition.\vspace*{-3pt}

\begin{proposition}[(Bad events are rare)]\label{bad-event} Suppose that
$p_0/2 \leq p < p_0$ and $l \leq\eps p/10$, and set $\delta:=
n_0^{-\kappa}$ for some sufficiently small fixed $\kappa> 0$.
Then:
\begin{longlist}[(b)]
\item[(a)] The events (\textup{i}), (\textup{iii}), (\textup{iv}), (\textup{v}) in Lemma \ref{backprop}
all fail with high probability.\vadjust{\goodbreak}
\item[(b)] There is a constant $C'$ such that all the coefficients of
the right singular vectors $u_j(M_{p,n})$ for $\eps p/2 \leq j \leq
(1-\eps/2) p$ are of magnitude at most $n^{-1/2} \log^{C'} n$ with
overwhelming probability. Conditioning $M_{p,n}$ to be a matrix with
this property, the events (\textup{ii}) and (\textup{vi}) occur with a conditional
probability of at most $2^{-\kappa m} + n^{-\kappa}$.
\item[(c)] Furthermore, there is a constant $C_2$ (depending on
$C',\kappa,C_1$) such that if $l \geq C_2$ and $M_{p,n}$ is
conditioned as in (\textup{b}), then (\textup{ii}) and (\textup{vi}) in fact occur with a
conditional probability of at most $2^{-\kappa m} \log^{-2C_1} n +
n^{-\kappa}$.\vspace*{-3pt}
\end{longlist}
\end{proposition}

But Proposition \ref{bad-event} can be proven by repeating the proof
of \cite{TVlocal1}, Proposition 53, with only cosmetic changes, the
only significant difference being that Theorem \ref{eigenconc} and
Theorem \ref{delocal} are applied instead of \cite{TVlocal1}, Theorem 60,
and \cite{TVlocal1}, Proposition 62, respectively.\vspace*{-3pt}

\section*{Acknowledgments}
We thank Horng-Tzer Yau for references, the anonymous referee for
helpful comments, and Jesse Geneson for a correction.  
\vspace*{-3pt}



\printaddresses

\end{document}